\RequirePackage{fix-cm}

\documentclass[smallextended]{svjour3}    

\smartqed  

\usepackage{latexsym}
\usepackage{amsmath,amssymb,amsfonts}

  \journalname{Positivity}

\begin{document}

\title{Universal function for a weighted space $L^1_{\mu}[0,1]$}

\author{Artsrun Sargsyan \and Martin Grigoryan}
         
            \institute{ A. Sargsyan\at
            YSU, Alex Manoogian 1, 0025, Yerevan, Armenia/CANDLE SRI, Acharyan 31, 0040, Yerevan, Armenia \\
              Tel.: +374-94-402155\\
              \email{asargsyan@ysu.am}
                              \and
              M. Grigoryan \at
              YSU, Alex Manoogian 1, 0025, Yerevan, Armenia \\
              Tel.: +374-77-456585\\
              \email{gmarting@ysu.am}          
}
\date{Received: date / Accepted: date}
\maketitle
\begin{abstract}
It is shown that there exist such a function $g\in L^1[0,1]$ and a weight function $0<\mu(x)\leq1$ that $g$ is universal for the weighted space $L^1_\mu[0,1]$ with respect to signs of its Fourier--Walsh coefficients.
\end{abstract}

\section{Introduction}
\label{intro}

Historically, the first type of a universal function was considered by G. Birkhoff \cite{Birkhoff1929} in 1929. He proved, that there exists an entire function $g(z)$, which is universal with respect to translations, i.e. for every entire function $f(z)$ and for each number $r>0$ there exists a growing sequence of natural numbers $\{n_k\}_{k=1}^\infty$, so that the sequence $\{g(z+n_k)\}_{k=1}^\infty$ uniformly converges to $f(z)$ on the disk $|z|\leq r$. In 1952 G. MacLane \cite{MacLane1952} proved a similar result for another type of universality, namely, there exists an entire function $g(z)$, which is universal with respect to derivatives, i.e. for every entire function $f(z)$ and for each number $r>0$ there exists a growing sequence of natural numbers $\{n_k\}_{k=1}^\infty$, so that the sequence $\{g^{(n_k)}(z)\}_{k=1}^\infty$ uniformly converges to $f(z)$ on the disk $|z|\leq r$. Further, in 1975 S. Voronin \cite{Voronin1975} proved the universality theorem for the Riemann zeta function $\zeta(s)$, which states that any nonvanishing analytic function can be approximated uniformly by certain purely imaginary shifts of the zeta function in the critical strip, namely, if $0<r<\frac{1}{4}$ and $g(s)$ is a nonvanishing continuous function on the disk $|s|\leq r$, that is analytic in the interior, then for any $\varepsilon>0$, there exists such a positive real number $\tau$ that
$$
\max_{|s|\leq r}\big|g(s)-\zeta(s+3/4+i\tau)\big|<\varepsilon .
$$
In 1987 K. Grosse--Erdman \cite{Grosse1987} proved the existence of infinitely differentiable function with universal Taylor expansion, namely, there exists such a function $g(x)\in C^{\infty}(\mathcal{R})$ with $g(0)=0$ that for every function $f(x)\in C(\mathcal{R})$ with $f(0)=0$ and for each number $r>0$ there exists a growing sequence of natural numbers $\{n_k\}_{k=1}^\infty$, so that the sequence
$$
S_{n_k}(g,0)=\sum_{m=1}^{n_k} {g^{(m)}(0)\over m!}x^{m}
$$
uniformly converges to $f(x)$ on $|x|\leq r$.

There are also many works devoted to the existence of universal series (in the common sense, with respect to rearrangements, partial series, signs of coefficients and etc.) in various classical orthogonal systems. The most general results were obtained by D. Menshov \cite{Menshov1964}, A. Talalyan \cite{Talalian19601}, P. Ulyanov \cite{Ul'yanov1972}  and their disciples (see \cite{Olevsky1968}--\cite{Episkoposian2006}).

The results presented in the current paper are an addition to this attractive area of mathematical research.
 
Let $|E|$ be the Lebesgue measure of a measurable set $E\subseteq [0,1]$, $\chi_E(x)$-- its characteristic function, $L^{p}[0,1]$ $(p>0)$ -- the class of all those measurable functions on $[0,1]$ that satisfy  the condition $\int_0^1|f(x)|^{p}dx<\infty$, $L^p_\mu[0,1]$ (weighted space) -- the class of all those measurable functions on $[0,1]$ that satisfy the condition $\int_0^1|f(x)|^p\mu (x)dx<\infty$, where $0<\mu(x)\leq1$ is a weight function, and $\{\varphi_k\}$ -- a complete orthonormal system in $L^2[0,1]$.

\paragraph{\bf Definition 1.}{\it
We say that a function $g\in L^{1}[0, 1]$  is universal for a class $L^{p}[0,1]$, $p>0$, with respect to signs of its Fourier coefficients (signs--subseries of its Fourier series) by the system $\{\varphi_k\}$, if for each function $f\in L^{p}[0,1]$ one can choose such numbers $\delta_{k}=\pm1$ ($\delta_{k}=0,\pm1$) that the series $\sum_{k=0}^{\infty}\delta_{k}c_{k}(g)\varphi_{k}(x)$, where $c_{k}(g)=\int
_{0}^{1}g(x)\varphi_{k}(x)dx$, converges to $f$ in $L^{p}[0,1]$ metric, i.e.
$$
\lim_{m\to \infty}\int_0^1\left|\sum_{k=0}^m\delta_{k}c_k(g)\varphi_k(x)-f(x)\right|^pdx=0.
$$
}

\paragraph{\bf Definition 2.} {\it Let $\mu(x)$ be a weight function defined on $[0,1]$. We say that a function $g\in L^1[0, 1]$  is universal for a weighted space $L^p_\mu [0,1]$ with respect to signs of its Fourier coeffcients (signs--subseries of its Fourier series) by the system $\{\varphi_k\}$, if for each function $f\in L^p_\mu [0,1]$ one can choose such numbers $\delta_{k}=\pm1$ ($\delta_{k}=0,\pm1$) that the series $\sum_{k=0}^{\infty}\delta_{k}c_{k}(g)\varphi_{k}(x)$ converges to $f$ in $L^p_\mu [0,1]$ metric, i.e.
$$
\lim_{m\to \infty}\int_0^1\left|\sum_{k=0}^m\delta_{k}c_k(g)\varphi_k(x)-f(x)\right|^p\mu (x)dx=0.
$$}

Let us recall the definition of the Walsh orthonormal system $\{W_n(x)\}_{n=0}^\infty$. Functions of the Walsh system are defined by means of Rademacher's functions 
$$
R_n(x)=\hbox{ sign}(\sin2^n\pi x ), \quad     x\in [0,1], \quad    n=1,2,\dots ,
$$
in the following way (see \cite{Golubov1987}): $W_0(x)\equiv 1$ and for $n\ge 1$
$$
 W_n(x)=\prod_{i=1}^p R_{k_i+1}(x),
$$
where $n=2^{k_1}+2^{k_2}+\dots+2^{k_p}\quad (k_1>k_2>\dots>k_p).$ 

\paragraph{\bf Remark.}{There does not exist a function $g\in L^1[0,1]$ which is universal for a certain class $L^p[0,1]$, $p\geq 1$, neither with respect to signs of its Fourier--Walsh coefficients nor with respect to signs--subseries of its Fourier--Walsh series.
}

Indeed, if such universal function existed then for the function $k_0c_{k_0}(g)W_{k_0}(x)$, where $k_0>1$ is any natural number with condition $c_{k_0}(g)\neq 0$, one could find such numbers $\delta_{k}=\pm1$ or $\delta_{k}=0,\pm1$ that
$$
\lim_{m\to \infty}\int_0^1\left|\sum_{k=0}^m\delta_kc_k(g)W_k(x) - k_0c_{k_0}(g)W_{k_0}(x)\right|^pdx=0,
$$
which simply leads to contradiction: $\delta_{k_0}=k_0>1$. 

It turns out, however, that the situation changes when considering weighted spaces $L^p_\mu[0, 1]$, $p\geq1$, or classes $L^p[0,1],\ p\in (0,1)$. For the latter case in \cite{GrigSar2016} the authors have constructed a universal function with respect to signs of Fourier--Walsh coefficients. For the case $L^p_\mu[0, 1]$, $p\geq1$ in \cite{GrigsSar2016} the authors could construct a universal function with respect to signs--subseries of Fourier--Walsh series. As a next step, the existence of a universal function for spaces $L^p_\mu[0, 1]$, $p\geq1$ with respect to signs of Fourier--Walsh coefficients was considered. We would like to note that the universality with respect to signs of Fourier--Walsh coefficients is much harder to acheive than the universality with respect to signs--subseries of Fourier--Walsh series. Therefore, it has been possible to construct such a function only for the $p=1$ case yet. The question is open for spaces $L^p_\mu[0, 1]$, $p>1$.

 The following theorem is true for the Walsh system:

\paragraph{\bf Theorem.}{\it There exist a function $g\in L^1[0,1]$ and a weight function $0<\mu(x)\leq1$, so that $g$ is universal for the weighted space $L^1_\mu[0,1]$ with respect to signs of its Fourier-Walsh coefficients.}

Moreover, it will be shown that the measure of the set on which $\mu(x)=1$ can be made arbitrarily close to 1, and
the function $g(x)$ can be choosen to have strictly decreasing Fourier--Walsh coefficients and converging to it by $L^1[0,1]$ norm Fourier--Walsh series.

Note that, considering generalities of many results obtained for the Walsh and trigonometric systems, an interest arises to find out whether the proved theorem is true for the trigonometric system.

\section{Main lemmas}
\label{Lem}
For the sake of simplicity the proof of the main theorem is divided into several steps which are given in the form of lemmas. Let us start from the known properties of the Walsh system $\{W_k(x)\}_{k=0}^\infty$. It is known (see \cite{Golubov1987}) that for each natural number $m$ we have
$$
\sum_{k=0}^{2^m-1}W_k(x)=
 \begin{cases}
 2^m, & \hbox  {when} \quad  x\in[0, 2^{-m}),
 \\
  0, &\hbox {when} \quad  x\in (2^{-m}, 1],
 \end{cases}\leqno(1.a)
$$
 and, consequently,
$$
\sum_{k=2^{m}}^{2^{m+1}-1}W_k(x)=
 \begin{cases}
 2^m, &\hbox {when} \quad  x\in [0, 2^{-m-1}),
 \\
 -2^m, & \hbox  {when} \quad  x\in(2^{-m-1}, 2^{-m}),
 \\
 0, &\hbox {when} \quad  x\in (2^{-m}, 1].
 \end{cases}\leqno(1.b)
$$

Obviously, for any natural number $M\in[2^m,2^{m+1})$ and numbers $\left\{a_k\right\}^{2^{m+1}-1}_{k=2^m}$ one gets
$$
\int_0^1\left|\sum_{k=2^{m}}^M a_k W_k(x)\right|dx\leq\left(\int_0^1\left| \sum_{k=2^{m}}^{2^{m+1}-1} a_k W_k(x)\right|^2dx\right)^{\frac{1}{2}}.\leqno(2)
$$
We will use also the following lemma \cite{Navasardyan1994}:

\paragraph{\bf Lemma 1.} {\it For each dyadic interval $\Delta=\left[\frac{l}{2^K},\frac{l+1}{2^K}\right]$, $l\in[0,2^K)$, and for every such natural number $M>K$ that $\frac{M-K}{2}$ is a whole number, there exists a polynomial in the Walsh system
$$
H(x)=\sum_{k=2^M}^{2^{M+1}-1}a_kW_k(x),
$$
so that

1) $|a_k|=2^{-\frac{M+K}{2}},\quad$ when  $\quad k\in[2^M,2^{M+1})$,

2) $H(x)=-1,\quad$ if $\quad x\in E_1,\ |E_1|=\frac{1}{2}|\Delta|$,

3) $H(x)=1,\quad$ if $\quad x\in E_2,\ |E_2|=\frac{1}{2}|\Delta|$,

4) $H(x)=0,\quad$ if $\quad x\not\in \Delta$.\newline
where $E_1$ and $E_2$ are finite unions of dyadic intervals.}

Now let us proceed to main lemmas of the paper:

\paragraph{\bf Lemma 2}{\it Let number $n_0\in\mathbb{N}$ and dyadic interval $\Delta=\left[\frac{l}{2^K},\frac{l+1}{2^K}\right]$, $l\in[0,2^K)$ be given. Then for any numbers $\varepsilon\in(0,1)$, $\gamma\neq 0$ and natural number $q $ there exist a measurable set $E_q\subset\Delta$  with measure $|E_q| =(1-2^{-q})|\Delta|$ and polynomials 
$$
P_q(x)=\sum_{k=2^{n_0}}^{{2^{n_q}}-1}a_kW_k(x)\quad\hbox{and}\quad H_q(x)=\sum_{k=2^{n_0}}^{{2^{n_q}}-1}\delta_ka_kW_k(x),\quad \delta_{k}=\pm1,
$$
in the Walsh system, so that 
$$
0<a_{k+1}\leq a_k<\varepsilon \quad \hbox{when} \quad k\in[2^{n_0}, 2^{n_q}-1),\leqno1)
$$ 
$$\int_{E_q}|\gamma\chi_{\Delta}(x)-H_q(x)|dx<\varepsilon\quad\hbox{and}\quad
\int_{[0,1]\setminus\Delta}|H_q(x)|dx<\varepsilon,\leqno2)
$$
$$
\max_{2^{n_0}\leq M < 2^{n_q}}\int_0^1\left| \sum_{k=2^{n_0}}^M \delta_ka_k W_k(x)\right|dx< 3|\gamma||\Delta|+\varepsilon,\leqno3)
$$
$$
\max_{2^{n_0}\leq M < 2^{n_q}}\int_0^1\left| \sum_{k=2^{n_0}}^Ma_k W_k(x)\right|dx<\varepsilon.\leqno4)
$$}

\paragraph{\bf Proof of lemma 2} The proof is performed by using the mathematical induction method with respect to the number $q$. We choose such a natural number $K_1>K$ that
$$
|\gamma|2^{-\frac{K_1+1}{2}}<\frac{\varepsilon}{2},\leqno(3)
$$
and present the interval $\Delta$ in the form of disjoint dyadic intervals' union
$$
\Delta=\bigcup_{i=1}^{N_1}\Delta_i^{(1)}
$$
with $\big|\Delta_i^{(1)}\big|=2^{-K_1-1},\ i\in[1,N_1].$
Obviously, $N_1=2^{K_1-K+1}$.

By denoting $K_{0}^{(1)}\equiv n_0-1$, for each natural number $i\in[1,N_1]$ we choose such a natural number $K_i^{(1)}>K_{i-1}^{(1)}\ \bigl(K_1^{(1)}>K_1\bigr)$ that the following conditions take place: 
\smallskip

a) $\frac{K_i^{(1)}-K_1-1}{2}$ is a whole number, 
\smallskip

b) $(K_i^{(1)}-K_{i-1}^{(1)})|\gamma|2^{-\frac{K_i^{(1)}+K_1+1}{2}}<\frac{\varepsilon}{4N_1},$
\smallskip

c) $2|\gamma|2^{-\frac{K_i^{(1)}+1}{2}}<\frac{\varepsilon}{2}.$
\smallskip

It immediately follows from (3) that
$$|\gamma|2^{-\frac{K_1^{(1)}+K_1+1}{2}}<\varepsilon.\leqno(4)$$

By successively applying lemma 1 for each interval $\Delta_i^{(1)}\ (i\in[1,N_1])$ and corresponding number $K_i^{(1)}$ we can find polynomials in the Walsh system 
$$
\widetilde H_i^{(1)}(x)=\sum_{k=2^{K_i^{(1)}}}^{2^{K_i^{(1)}+1}-1}\tilde a_kW_k(x),\quad i\in[1,N_1],\leqno(5)
$$
so that
$$
|\tilde a_k| =|\gamma|2^{-\frac{K_i^{(1)} + K_1+1}{2}}\quad \hbox{when}\quad k\in\bigl[2^{K_i^{(1)}}, 2^{K_i^{(1)}+1}\bigr),\leqno(6)
$$
$$
\widetilde H_i^{(1)}(x)=
\begin{cases}
-\gamma, & \hbox{for} \quad x\in \widetilde{E_i}^{(1)}\subset \Delta_i^{(1)},\quad \big|\widetilde{E_i}^{(1)}\big|=\frac{\big|\Delta_i^{(1)}\big|}{2},
\\
\gamma, & \hbox{for} \quad x\in \widetilde{\widetilde{E_i}}^{(1)}\subset \Delta_i^{(1)},\quad \big|\widetilde{\widetilde{E_i}}^{(1)}\big|=\frac{\big|\Delta_i^{(1)}\big|}{2},
\\
0, & \hbox{for} \quad x\notin \Delta_i^{(1)}.
\end{cases}\leqno(7)
$$
Hence, by denoting
$$
\widetilde H_1(x)=\sum_{i=1}^{N_1}\widetilde H_i^{(1)}(x),\leqno(8)
$$
we get
$$
\widetilde H_1(x)=
\begin{cases}
-\gamma, & \hbox{for} \quad x\in \widetilde E_1\subset \Delta,\quad  \big|\widetilde E_1\big|= \frac{|\Delta|}{2},
\\
\gamma, & \hbox{for} \quad x\in \Delta\setminus \widetilde E_1,
\\
0, & \hbox{for} \quad x\notin \Delta.
\end{cases}\leqno(9)
$$

As the polynomial $\widetilde H_i^{(1)}(x)$ is a linear combination of Walsh functions from $K_i^{(1)}$ group, it is clear, that the set $\widetilde E_1$ can be presented as a union of certain $N_2$ amount of disjoint dyadic intervals
$\Delta_i^{(2)}$ with measure $\big|\Delta_i^{(2)}\big|=2^{-K_{N_1}^{(1)} -1},\ i\in[1,N_2]$:
$$
\widetilde E_1=\bigcup_{i=1}^{N_2}\Delta_i^{(2)}.
$$

After making the following definitions
$$
E_1 = \Delta\setminus \widetilde E_1\leqno(10)
$$
and 
$$
\begin{cases}
\tilde a_k = |\gamma|2^{-\frac{K_i^{(1)}+K_1+1}{2}},& \hbox{when}\quad k\in \bigl[2^{K_{i-1}^{(1)}+1}, 2^{K_i^{(1)}}\bigl),\quad i\in[1,N_1],
\\
a_k=|\tilde a_k|,\quad \delta_k=\frac{\tilde a_k}{|\tilde a_k|},& \hbox{when}\quad k\in\bigl[2^{n_0},2^{K_{N_1}^{(1)}+1}\bigr),
\end{cases}\leqno(11)
$$
let us verify that the set $E_1$ and polynomials
$$
P_1(x)=\sum_{k=2^{n_0}}^{{2^{K_{N_1}^{(1)}+1}}-1}a_kW_k(x),\quad H_1(x)=\sum_{k=2^{n_0}}^{{2^{K_{N_1}^{(1)}+1}}-1}\delta_ka_kW_k(x)
$$
satisfy all lemma 2 statements for $q=1$. Indeed, from (9) and (10) it immediately follows that $| E_1|=(1-2^{-1})|\Delta|$. The statement 1) follows from (4), (6), (11) and from monotonicity of numbers $K^{(1)}_i\ (i\in[1,N_1])$. For the proof of statements 2) and 3) we present the polynomial $H_1(x)$ in the form of
$$
H_1(x)=\widetilde H_1(x)+\sum_{i=1}^{N_1} \sum_{k=2^{K_{i-1}^{(1)}+1}}^{2^{K_{i}^{(1)}}-1} a_kW_k(x)=
\leqno(12)
$$
$$
=\widetilde H_1(x)+\sum_{i=1}^{N_1}\sum_{j=1}^{K_{i}^{(1)}-K_{i-1}^{(1)}-1} Q_{i,j}^{(1)}(x),
$$
where
$$
Q_{i,j}^{(1)}(x)=\sum_{k=2^{K_{i-1}^{(1)}+j}}^{2^{K_{i-1}^{(1)}+j+1}-1} a_kW_k(x)\leqno(13)
$$
$\bigl(\delta_k =1\ \hbox{when}\ k\in \bigl[2^{K_{i-1}^{(1)}+1}, 2^{K_i^{(1)}}\bigr), i\in[1,N_1]\bigr)$. 

Since all coefficients $a_k$ are equal when  $k\in \bigl[2^{K_{i-1}^{(1)}+1}, 2^{K_i^{(1)}}\bigl),\ i\in[1,N_1]$ (see (11)) then considering (1.b), (9)--(13) and b) condition for numbers $K_{i}^{(1)}\ (i\in[1,N_1])$ we obtain
$$
\int_{E_1}\left|\gamma\chi_\Delta(x)-H_1(x)\right|dx\leq \int_{E_1}\left|\gamma\chi_\Delta(x)-\widetilde H_1(x)\right|dx +
$$ 
$$
+\sum_{i=1}^{N_1}\sum_{j=1}^{K_{i}^{(1)}-K_{i-1}^{(1)}-1}\int_0^1 \big|Q_{i,j}^{(1)}(x)\big| dx
=\sum_{i=1}^{N_1}\bigl(K_{i}^{(1)}-K_{i-1}^{(1)}-1\bigr)|\gamma|2^{-\frac{K_i^{(1)}+K_1+1}{2}}<\varepsilon
$$
and
$$
\int_{[0,1]\setminus\Delta}\left|H_1(x)\right|dx\leq\sum_{i=1}^{N_1}\sum_{j=1}^{K_{i}^{(1)}-K_{i-1}^{(1)}-1}\int_0^1 \big|Q_{i,j}^{(1)}(x)\big| dx
<\varepsilon.
$$

To prove statements 3) and 4) we present the natural number $M\in \bigl[2^{n_0}, 2^{K_{N_1}^{(1)}+1}\bigr)$ in the form of $M = 2^{\tilde n} +s,\ s\in[0, 2^{\tilde n})$, where $\tilde n\in \bigl(K_{m-1}^{(1)}, K_{m}^{(1)}\bigr]$ for some $m\in[1,N_1]$. From (8), (12) and (13) it follows that
$$
\int_0^1\left| \sum_{k=2^{n_0}}^M \delta_ka_k W_k(x)\right|dx\leq \int_0^1\left|\sum_{i=1}^{m-1}\widetilde H_i^{(1)}(x)\right|dx+
$$
$$
+\sum_{i=1}^{m-1}\sum_{j=1}^{K_{i}^{(1)}-K_{i-1}^{(1)}-1}\int_0^1 \big|Q_{i,j}^{(1)}(x)\big| dx+
$$
$$
+\sum_{j=1}^{\tilde n-K_{m-1}^{(1)}-1}\int_0^1 \big|Q_{m-1,j}^{(1)}(x)\big|dx+\int_0^1 \left| \sum_{k=2^{\tilde n}}^{2^{\tilde n}+s}\delta_k a_kW_k(x)\right| dx.
$$
By using (1.b)--(3), (5)--(8) and (11) we get
$$
\int_0^1\left|\sum_{i=1}^{m-1}\widetilde H_i^{(1)}(x)\right|dx + \int_0^1 \left| \sum_{k=2^{\tilde n}}^{2^{\tilde n}+s}\delta_k a_kW_k(x)\right|dx\leq
$$
$$
\leq\sum_{i=1}^{m-1}\int_0^1\left|\widetilde H_i^{(1)}(x)\right|dx+\left(\int_0^1 \left| \sum_{k=2^{\tilde n}}^{2^{\tilde n}+s}\delta_k a_kW_k(x)\right|^2 dx\right)^{\frac{1}{2}}\leq
$$
$$
\leq\sum_{i=1}^{N_1}\int_0^1\left|\widetilde H_i^{(1)}(x)\right|dx+\left(\int_0^1\left| \sum_{k=2^{\tilde n}}^{2^{\tilde n +1}-1} \delta_ka_k W_k(x)\right|^2dx\right)^{\frac{1}{2}}=
$$
$$
=|\gamma||\Delta|+|\gamma|2^{-\frac{K_m^{(1)}+K_1+1}{2}}2^{\frac{\tilde n}{2}}<3|\gamma||\Delta|+\frac{\varepsilon}{2},
$$
thus, taking 
$$
\sum_{i=1}^{m-1}\sum_{j=1}^{K_{i}^{(1)}-K_{i-1}^{(1)}-1}\int_0^1 \big|Q_{i,j}^{(1)}(x)\big| dx+\sum_{j=1}^{\tilde n-K_{m-1}^{(1)}-1}\int_0^1 \big|Q_{m-1,j}^{(1)}(x)\big|dx\leq
$$
$$
\leq\sum_{i=1}^{N_1}\sum_{j=1}^{K_{i}^{(1)}-K_{i-1}^{(1)}-1}\int_0^1 \big|Q_{i,j}^{(1)}(x)\big| dx
=\sum_{i=1}^{N_1}\bigl(K_{i}^{(1)}-K_{i-1}^{(1)}-1\bigr)|\gamma|2^{-\frac{K_i^{(1)}+K_1+1}{2}}<\frac{\varepsilon}{4}
$$
into account we verify the validity of the statement 3).

Further, for each natural number $n\in \bigl[n_0, K_{N_1}^{(1)}\bigr]$ we denote $b_n=a_k,\ k\in[2^n, 2^{n+1})$ (coefficients $a_k$ of Walsh functions from n--th group are equal in $H_1(x)$). It follows from (1.b)--(3), (6), (11) and b) condition for numbers $K_{i}^{(1)}\ (i\in[1,N_1])$ that 
$$
\sum_{n={n_0}}^{K_{N_1}^{(1)}}b_n=\sum_{i=1}^{N_1} \sum_{n={K_{i-1}^{(1)}+1}}^{{K_{i}^{(1)}}} b_n=\sum_{i=1}^{N_1}\bigl(K_{i}^{(1)}-K_{i-1}^{(1)}\bigr)|\gamma|2^{-\frac{K_{i}^{(1)}+K_1+1}{2}}<\frac{\varepsilon}{4}
$$
and
$$
\int_0^1\left| \sum_{k=2^{n_0}}^Ma_k W_k(x)\right|dx\leq\sum_{n={n_0}}^{\tilde n-1}b_n + \int_0^1\left| \sum_{k=2^{\tilde n}}^{2^{\tilde n}+s} a_k W_k(x)\right|dx\leq
$$
$$
\leq\sum_{n={n_0}}^{K_{N_1}^{(1)}}b_n+
\left(\int_0^1\left|\sum_{k=2^{\tilde n}}^{2^{\tilde n+1}-1} b_{\tilde n} W_k(x)\right|^2\right)^{\frac{1}{2}}<\frac{\varepsilon}{4}+|\gamma|2^{-\frac{K_{m}^{(1)}+K_1+1}{2}}2^{\frac{\tilde n}{2}}<\varepsilon,
$$
which proves the statement 4).

Assume that for $q>1$ natural numbers 
$$
K^{(1)}_1<\dots<K^{(1)}_{N_{1}}<\dots<K^{(q-1)}_1<\dots<K^{(q-1)}_{N_{q-1}},$$
sets 
$$
\widetilde E_{q-1}\subset \Delta\quad \hbox{and}\quad E_{q-1}=\Delta\setminus \widetilde E_{q-1}
$$
and polynomials
$$
P_{q-1}(x)=\sum_{k=2^{n_0}}^{{2^{K_{N_{q-1}}^{(q-1)}+1}}-1}a_kW_k(x),\quad H_{q-1}(x)=\sum_{k=2^{n_0}}^{{2^{K_{N_{q-1}}^{(q-1)}+1}}-1}\delta_ka_kW_k(x),\ \delta_k=\pm 1
$$ 
are already chosen to satisfy the following conditions:
\smallskip

$a^\prime$) $\frac{K_i^{(\nu)}-K_{N_{\nu-1}}^{(\nu-1)}-1}{2}$ is a whole number
$\bigr(K_{N_{0}}^{(0)}\equiv K_1\bigl),$ 
\smallskip

$b^\prime$) $\bigl(K_i^{(\nu)}-K_{i-1}^{(\nu)}\bigr)2^{\nu-1}|\gamma|2^{-\frac{K_i^{(\nu)}+K_{N_{\nu-1}}^{(\nu-1)}+1}{2}}<\frac{\varepsilon}{2^{\nu +1}N_{\nu}},$
\smallskip

$c^\prime$) $2^{\nu}|\gamma|2^{-\frac{K_i^{(\nu)}+1}{2}}<\frac{\varepsilon}{2},$

$$
a_k = 2^{\nu-1}|\gamma|2^{-\frac{K_i^{(\nu)}+K_{N_{\nu-1}}^{(\nu-1)}+1}{2}}\quad \hbox{for}\quad k\in \bigl[2^{K_{i-1}^{(\nu)}+1}, 2^{K_i^{(\nu)}+1}\bigl),\leqno(14)
$$
$$
K_0^{(\nu)}\equiv 
\begin{cases}
K_{N_{\nu-1}}^{(\nu-1)}, & \hbox{if} \quad \nu>1,
\\
n_0-1, & \hbox{if} \quad \nu=1,
\end{cases}
$$
for any natural numbers $i\in[1,N_{\nu}]$ and $\nu\in [1, q-1]$. Besides,
$$
\sum_{n={n_0}}^{K_{N_{q-1}}^{(q-1)}} b_n<\sum_{k=1}^{q-1}\frac{\varepsilon}{2^{k+1}},\quad\hbox{where}\quad b_n\equiv a_k,\ k\in[2^n, 2^{n+1}),\leqno(15)
$$
$$
H_{q-1}(x)=\widetilde H_{q-1}(x)+\sum_{\nu=1}^{q-1}\sum_{i=1}^{N_\nu} \sum_{k=2^{K_{i-1}^{(\nu)}+1}}^{2^{K_{i}^{(\nu)}}-1} a_kW_k(x),\leqno(16)
$$
$$
\widetilde H_{q-1}(x)=
\begin{cases}
-(2^{q-1}-1)\gamma,& \hbox{for} \quad x\in \widetilde E_{q-1},
\\
\gamma, & \hbox{for} \quad x\in E_{q-1},
\\
0,& \hbox{for} \quad x\notin \Delta,
\end{cases}\leqno(17)
$$
$$
\big|\widetilde E_{q-1}\big|=2^{-q+1}|\Delta|\quad \hbox{and}\quad \big|E_{q-1}\big|= \bigr(1-2^{-q+1}\bigr)|\Delta|\leqno(18)
$$
and the set  $\widetilde E_{q-1}$ can be presented as a union of certain $N_q$ amount of disjoint dyadic intervals $\Delta_i^{(q)}$ with measure $\big|\Delta_i^{(q)}\big|=2^{-K^{(q-1)}_{N_{q-1}}-1},\ i\in[1,N_{q}]$:
$$
\widetilde E_{q-1}=\bigcup_{i=1}^{N_{q}}\Delta_i^{(q)}.
$$

For each natural number $i\in[1,N_{q}]$ we choose such a natural number $K_i^{(q)}>K_{i-1}^{(q)}$ $\bigl(K_{0}^{(q)}\equiv K_{N_{q-1}}^{(q-1)}\bigr)$ that the following conditions hold: 
\smallskip

$a^{\prime\prime})$ $\frac{K_i^{(q)}-K_{N_{q-1}}^{(q-1)}-1}{2}$ is a whole number, 
\smallskip

$b^{\prime\prime})$ $\bigl(K_i^{(q)}-K_{i-1}^{(q)}\bigr)2^{q-1}|\gamma|2^{-\frac{K_i^{(q)}+K_{N_{q-1}}^{(q-1)}+1}{2}}<\frac{\varepsilon}{2^{q+1}N_{q}},$
\smallskip

$c^{\prime\prime})$ $2^{q}|\gamma|2^{-\frac{K_i^{(q)}+1}{2}}<\frac{\varepsilon}{2}.$
\smallskip

By a successive application of lemma 1 for each interval $\Delta_i^{(q)}\subset \widetilde E_{q-1}\ (i\in[1,N_{q}])$ and corresponding number $K_i^{(q)}$ we can find polynomials in the Walsh system 
$$
\widetilde H_i^{(q)}(x)=\sum_{k=2^{K_i^{(q)}}}^{2^{K_i^{(q)}+1}-1}\tilde a_kW_k(x),\quad i\in[1,N_{q}],\leqno(19)
$$
so that
$$
|\tilde a_k|= 2^{q-1}|\gamma|2^{-\frac{K_i^{(q)} + K_{N_{q-1}}^{(q-1)}+1}{2}}\quad  \hbox{when} \quad k\in \bigl[2^{K_i^{(q)}}, 2^{K_i^{(q)}+1}\bigr),\leqno(20)
$$
$$
\widetilde H_i^{(q)}(x)=
\begin{cases}
-2^{q-1}\gamma, & \hbox{for} \quad x\in \widetilde{E_i}^{(q)}\subset \Delta_i^{(q)},\quad \big|\widetilde{E_i}^{(q)}\big|=\frac{\big|\Delta_i^{(q)}\big|}{2},
\\
2^{q-1}\gamma, & \hbox{for} \quad x\in \widetilde{\widetilde{E_i}}^{(q)}\subset \Delta_i^{(q)},\quad \big|\widetilde{\widetilde{E_i}}^{(q)}\big| =\frac{\big|\Delta_i^{(q)}\big|}{2},
\\
0, & \hbox{for} \quad x\notin \Delta_i^{(q)}.
\end{cases}\leqno(21)
$$
Hence, by denoting
$$
\widetilde H_q(x)=\widetilde H_{q-1}(x)+\sum_{i=1}^{N_{q}}\widetilde H_i^{(q)}(x)\leqno(22)
$$
and taking (17) and (18) into account we obtain
$$
\widetilde H_q(x)=
\begin{cases}
-(2^q-1)\gamma, & \hbox{for} \quad x\in \widetilde E_q\subset \widetilde E_{q-1},\quad  \big|\widetilde E_q\big|=\frac{|\Delta|}{2^q},
\\
\gamma, & \hbox{for} \quad x\in \Delta\setminus \widetilde E_q,
\\
0, & \hbox{for} \quad x\notin \Delta.
\end{cases}\leqno(23)
$$

After defining
$$
E_q=\Delta\setminus \widetilde E_q\leqno(24)
$$
and
$$
\begin{cases}
\tilde a_k = 2^{q-1}|\gamma|2^{-\frac{K_i^{(q)}+K_{N_{q-1}}^{(q-1)}+1}{2}},& \hbox{for}\   k\in \bigl[2^{K_{i-1}^{(q)}+1}, 2^{K_i^{(q)}}\bigr),\ i\in [1, N_{q}],
\\
a_k=|\tilde a_k|,\quad \delta_k=\frac{\tilde a_k}{|\tilde a_k|} ,&\hbox{for}\  k\in\bigl[2^{n_0},2^{K_{N_q}^{(q)}+1}\bigr),
\end{cases}\leqno(25)
$$
let us verify that the set $E_q$ and polynomials 
$$
P_q(x)=\sum_{k=2^{n_0}}^{{2^{n_q}}-1}a_kW_k(x),\quad H_q(x)=\sum_{k=2^{n_0}}^{{2^{n_q}}-1}\delta_ka_kW_k(x),
$$
where $n_q\equiv K_{N_{q}}^{(q)}+1$, satisfy all lemma 2 statements. Indeed, from (23) and (24) it follows that $\big|E_q\big|= (1-2^{-q})|\Delta|$. The statement 1) follows from (4), (14), (20), (25) and from monotonicity of numbers $K^{(\nu)}_i\ (i\in[1,N_{\nu}],\ \nu\in[1,q])$. For the proof of statements 2) and 3) we present the polynomial $H_q(x)$ in the form of
$$
H_q(x)=\widetilde H_q(x)+\sum_{\nu=1}^{q}\sum_{i=1}^{N_\nu} \sum_{k=2^{K_{i-1}^{(\nu)}+1}}^{2^{K_{i}^{(\nu)}}-1} a_kW_k(x)=\leqno(26)
$$
$$
=\widetilde H_q(x)+\sum_{\nu=1}^{q}\sum_{i=1}^{N_\nu}\sum_{j=1}^{K_{i}^{(\nu)}-K_{i-1}^{(\nu)}-1} Q_{i,j}^{(\nu)}(x),
$$
where  
$$
Q_{i,j}^{(\nu)}(x)=\sum_{k=2^{K_{i-1}^{(\nu)}+j}}^{2^{K_{i-1}^{(\nu)}+j+1}-1} a_kW_k(x).\leqno(27)
$$

Considering relations (1.b), (14), (23)--(26) and $b^\prime)$, $b^{\prime\prime})$ conditions for numbers $K_i^{(\nu)}\ \left(i\in[1,N_\nu],\ \nu\in[1,q]\right)$ we have
$$
\int_{E_q}\left|\gamma\chi_\Delta(x)-H_q(x)\right|dx\leq \int_{E_q}\left|\gamma\chi_\Delta(x)-\widetilde H_q(x)\right|dx +
$$ 
$$
+\sum_{\nu=1}^q\sum_{i=1}^{N_\nu}\sum_{j=1}^{K_{i}^{(\nu)}-K_{i-1}^{(\nu)}-1}\int_0^1 \big|Q_{i,j}^{(\nu)}(x)\big| dx=
$$
$$=\sum_{\nu=1}^{q}\sum_{i=1}^{N_{\nu}} \bigl(K_{i}^{(\nu)}-K_{i-1}^{(\nu)}-1\bigr)2^{\nu-1}|\gamma|2^{-\frac{K_i^{(\nu)}+K_{N_{\nu-1}}^{(\nu-1)}+1}{2}}<\varepsilon
$$
and
$$
\int_{[0,1]\setminus\Delta}\left|H_q(x)\right|dx\leq\sum_{\nu=1}^{q}\sum_{i=1}^{N_\nu}\sum_{j=1}^{K_{i}^{(\nu)}-K_{i-1}^{(\nu)}-1}\int_0^1 \big|Q_{i,j}^{(\nu)}(x)\big| dx
<\varepsilon.
$$

To prove statements 3) and 4) we present the natural number $M\in [2^{n_0}, 2^{n_q})$ in the form of $M=2^{\tilde n} + s,\ s\in[0, 2^{\tilde n})$. Let us consider only the case when $\tilde n\in\bigl(K_{N_{q-1}}^{(q-1)}, K_{N_{q}}^{(q)}\bigr]$, since all other cases were under consideration in previous steps of induction. Let $\tilde n\in\bigl(K_{m-1}^{(q)}, K_{m}^{(q)}\bigr]$ for some $m\in[1,N_q]$. From (16), (22), (26) and (27) it follows that
$$
\int_0^1\left| \sum_{k=2^{n_0}}^M \delta_ka_k W_k(x)\right|dx\leq \int_0^1\left|\widetilde H_{q-1}(x)+\sum_{i=1}^{m-1}\widetilde H_i^{(q)}(x)\right|dx +
$$
$$
+\sum_{\nu=1}^{q-1}\sum_{i=1}^{N_{\nu}} \sum_{j=1}^{K_{i}^{(\nu)}-K_{i-1}^{(\nu)}-1}\int_0^1 \left|Q_{i,j}^{(\nu)}(x)\right| dx
+\sum_{i=1}^{m-1}\sum_{j=1}^{K_{i}^{(q)}-K_{i-1}^{(q)}-1}\int_0^1 \left|Q_{i,j}^{(q)}(x)\right| dx+
$$
$$
+\sum_{j=1}^{\tilde n-K_{m-1}^{(q)}-1}\int_0^1 \left|Q_{m-1,j}^{(q)}(x)\right| dx+\int_0^1 \left| \sum_{k=2^{\tilde n}}^{2^{\tilde n}+s}\delta_k a_kW_k(x)\right| dx.
$$
By analyzing relations (1.b), (2), (17)--(21), (25) and $b^\prime)$ $b^{\prime\prime})$, $c^\prime)$ conditions for numbers $K_{i}^{(\nu)}\ \left(i\in[1,N_\nu],\ \nu\in[1,q]\right)$  we get
$$
\int_0^1\left|\widetilde H_{q-1}(x)+\sum_{i=1}^{m-1}\widetilde H_i^{(q)}(x)\right|dx+ \int_0^1 \left| \sum_{k=2^{\tilde n}}^{2^{\tilde n}+s}\delta_k a_kW_k(x)\right| dx\leq
$$
$$
\leq \int_0^1|\widetilde H_{q-1}(x)|dx+\sum_{i=1}^{m-1}\int_0^1|\widetilde H_i^{(q)}(x)|dx+\left(\int_0^1 \left| \sum_{k=2^{\tilde n}}^{2^{\tilde n}+s}\delta_k a_kW_k(x)\right|^2 dx\right)^{\frac{1}{2}}\leq
$$
$$
\leq\int_0^1|\widetilde H_{q-1}(x)|dx+\sum_{i=1}^{N_q}\int_0^1|\widetilde H_i^{(q)}(x)|dx+\left(\int_0^1\left| \sum_{k=2^{\tilde n}}^{2^{\tilde n +1}-1} \delta_ka_k W_k(x)\right|^2dx\right)^{\frac{1}{2}}<
$$
$$
<2|\gamma||\Delta|+|\gamma||\Delta|+2^{q-1}|\gamma|2^{-\frac{K_m^{(q)}+K_{N_{q-1}}^{(q-1)}+1}{2}}2^{\frac{\tilde n}{2}}<3|\gamma||\Delta|+\frac{\varepsilon}{2}
$$
and 
$$
\sum_{\nu=1}^{q-1}\sum_{i=1}^{N_{\nu}} \sum_{j=1}^{K_{i}^{(\nu)}-K_{i-1}^{(\nu)}-1}\int_0^1 \left|Q_{i,j}^{(\nu)}(x)\right| dx
+\sum_{i=1}^{m-1}\sum_{j=1}^{K_{i}^{(q)}-K_{i-1}^{(q)}-1}\int_0^1 \left|Q_{i,j}^{(q)}(x)\right| dx+
$$
$$
+\sum_{j=1}^{\tilde n-K_{m-1}^{(q)}-1}\int_0^1 \left|Q_{m-1,j}^{(q)}(x)\right| dx
\leq\sum_{\nu=1}^{q}\sum_{i=1}^{N_{\nu}} \sum_{j=1}^{K_{i}^{(\nu)}-K_{i-1}^{(\nu)}-1}\int_0^1 \left|Q_{i,j}^{(\nu)}(x)\right| dx=
$$
$$
=\sum_{\nu=1}^{q}\sum_{i=1}^{N_{\nu}} \bigl(K_{i}^{(\nu)}-K_{i-1}^{(\nu)}-1\bigr)2^{\nu-1}|\gamma|2^{-\frac{K_i^{(\nu)}+K_{N_{\nu-1}}^{(\nu-1)}+1}{2}}<\frac{\varepsilon}{2},
$$
which proves the statement 3).

Further, for each natural number $n\in\bigl[n_0,K_{N_{q}}^{(q)}\bigr]$ we denote
$$
b_n\equiv a_k,\quad\hbox{when}\quad k\in[2^n, 2^{n+1}),
$$
and use (1.b), (2), (15), (20), (25), $c^\prime)$ condition for number $K_{N_{q-1}}^{(q-1)}$ and $b^{\prime\prime})$ condition for numbers $K_{i}^{(q)} (i\in[1,N_q])$ to obtain
$$
\int_0^1\left| \sum_{k=2^{n_0}}^Ma_k W_k(x)\right|dx\leq\sum_{n={n_0}}^{\tilde n-1}b_n + \int_0^1\left| \sum_{k=2^{\tilde n}}^{2^{\tilde n}+s} a_k W_k(x)\right|dx\leq
$$
$$
\leq\sum_{n={n_0}}^{K_{N_{q-1}}^{(q-1)}} b_n +\sum_{i=1}^{N_{q}}\sum_{n={K_{i-1}^{(q)}+1}}^{{K_{i}^{(q)}}} b_n+
\left(\int_0^1\left|\sum_{k=2^{\tilde n}}^{2^{\tilde n+1}-1} b_{\tilde n} W_k(x)\right|^2\right)^{\frac{1}{2}}\leq
$$
$$
\leq\sum_{k=1}^{q-1}\frac{\varepsilon}{2^{k+1}}+\sum_{i=1}^{N_{q}}\bigl(K_{i}^{(q)}-K_{i-1}^{(q)}\bigr)2^{q-1}|\gamma|2^{-\frac{K_{i}^{(q)}+K_{N_{q-1}}^{(q-1)}+1}{2}}+
$$
$$
+2^{q-1}|\gamma|2^{-\frac{K_{m}^{(q)}+K_{N_{q-1}}^{(q-1)}+1}{2}}2^{\frac{\tilde n}{2}}<\varepsilon.
$$

Lemma 2 is proved.

\paragraph{\bf Lemma 3.}{\it Let $n_0\in\mathbb{N},\ \varepsilon\in(0,1)$ and $f(x)=\sum_{j=1}^{\nu_0} \gamma_j\chi_{\Delta_j}(x)$ be such a step function that $\gamma_j\neq 0$ and $\{\Delta_j\}_{j=1}^{\nu_0}$ are disjoint dyadic intervals with 
$\sum_{j=1}^{\nu_0}|\Delta_j|=1$. Then one can find a measurable set $E_\varepsilon\subset[0,1]$ with measure $|E_\varepsilon| > 1-\varepsilon$ and  polynomials 
$$
P(x)=\sum_{k=2^{n_{0}}}^{2^{n}-1}a_{k}W_{k}(x)\quad\hbox{and}\quad H(x)=\sum_{k=2^{n_{0}}}^{2^{n}-1}\delta_{k}a_{k}W_{k}(x),\quad\delta_{k}=\pm1,
$$
in the Walsh system, which satisfy the following conditions:
$$
0<a_{k+1}<a_{k}<\varepsilon ,\quad k\in\mathit{[\mathit{2}^{n_{0}},2^{n}-1)},\leqno1)
$$
$$
\int_{E_\varepsilon}\left|f(x)-H(x)\right|dx<\varepsilon, \leqno2)
$$
$$
\max_{2^{n_{0}}\leq M<2^{n}}\int_{e}\left|\sum_{k=2^{n_{0}}}^{M}\delta_{k}
a_{k}W_{k}(x)\right|dx<\int_e|f(x)|dx+\varepsilon\leqno3)
$$
for any measurable set $e\subseteq E_\varepsilon$,
$$
\max_{2^{n_{0}}\leq M<2^{n}}\int_{0}^{1}\left|\sum_{k=2^{n_{0}}}^{M}a_{k}
W_{k}(x)\right|dx<\varepsilon.\leqno4)
$$}

\paragraph{\bf Proof of lemma 3}
Without loss of generality we can assume that intervals $\Delta_j\ (j\in[1,\nu_0])$ have the same length and are small enough to privide the condition
$$
\max_{j\in[1,\nu_0]}3|\gamma_j||\Delta_j|<\frac{\varepsilon}{8}.\leqno(28)
$$

We choose a natural number $q>\log_2\frac{1}{\varepsilon}$ and apply  lemma 2 for each interval $\Delta_j$ to get sets $E_q^{(j)}\subset \Delta_j$ with measure
$$
\big|E_q^{(j)}\big|= (1-2^{-q})|\Delta_j|>(1-\varepsilon)|\Delta_j|\leqno(29)
$$
and polynomials
$$
\widetilde P_q^{(j)}(x)=\sum_{k=2^{n_{j-1}}}^{{2^{n_j}}-1}\tilde a_k^{(j)}W_k(x),\quad \widetilde H_q^{(j)}(x)=\sum_{k=2^{n_{j-1}}}^{{2^{n_j}}-1}\delta_k^{(j)}\tilde a_k^{(j)}W_k(x),\ \delta_k^{(j)}=\pm 1,
$$ 
in the Walsh system, which satisfy the following conditions:
$$
\begin{cases}
0<\tilde a_{k+1}^{(1)}\leq\tilde a_k^{(1)}<\frac{\varepsilon}{2}, &\hbox{for}\ k\in[2^{n_{0}}, 2^{n_1}-1),
\\
0<\tilde a_{k+1}^{(j)}\leq\tilde a_k^{(j)}<\tilde a_{{2^{n_{j-1}}}-1}^{(j-1)}, &\hbox{for}\ k\in[2^{n_{j-1}}, 2^{n_j}-1),\  j\in[2,\nu_0],
\end{cases}\leqno(30)
$$ 
$$
\int_{E_q^{(j)}}|\gamma_j\chi_{\Delta_j}(x)-\widetilde H_q^{(j)}(x)|dx<\frac{\varepsilon}{8\nu_0}\quad\hbox{and}\quad\int_{[0,1]\setminus\Delta_j}|\widetilde H_q^{(j)}(x)|dx<\frac{\varepsilon}{8\nu_0^2},\leqno(31)
$$
$$
\max_{2^{n_{j-1}}\leq M < 2^{n_j}}\int_0^1\left| \sum_{k=2^{n_{j-1}}}^M \delta_k^{(j)}\tilde a_k^{(j)} W_k(x)\right|dx<\leqno(32)
$$
$$
< 3|\gamma_j||\Delta_j|+\frac{\varepsilon}{8}<\frac{\varepsilon}{4}\quad(\hbox{see (28)})
$$
and
$$
\max_{2^{n_{j-1}}\leq M < 2^{n_j}}\int_0^1\left| \sum_{k=2^{n_{j-1}}}^M\tilde a_k^{(j)} W_k(x)\right|dx<\frac{\varepsilon}{2^{j+1}}.\leqno(33)
$$

We define a set $E_\varepsilon$ and polynomials $\widetilde P(x)$, $\widetilde H(x)$ in the following way:
$$
E_\varepsilon = \bigcup_{j=1}^{\nu_0}E_q^{(j)},\leqno(34)
$$
$$
\widetilde P(x) =\sum_{j=1}^{\nu_0}\widetilde P_q^{(j)}(x)=\sum_{k=2^{n_0}}^{2^{n_{\nu_0}-1}}\tilde a_kW_k(x),
$$
$$
\widetilde H(x) =\sum_{j=1}^{\nu_0}\widetilde H_q^{(j)}(x)=\sum_{k=2^{n_0}}^{2^{n_{\nu_0}-1}}\delta_k\tilde a_kW_k(x),
$$
where $\tilde a_k=\tilde a_k^{(j)}$ and $\delta_k=\delta_k^{(j)}$, when $k\in[2^{n_{j-1}},2^{n_j})$. 

It follows from (29), (30) and (34) that
$$
|E_\varepsilon| >1-\varepsilon,
$$
$$
0<\tilde a_{k+1}\leq \tilde a_k<\frac{\varepsilon}{2},\quad\hbox{when}\quad k\in[2^{n_0}, 2^{n_{\nu_0}}-1).
$$

Considering (31) and (34), we have
$$
\int_{E_\varepsilon}\left|f(x)-\widetilde H(x)\right|dx=\sum_{j=1}^{\nu_0}\int_{E_q^{(j)}}\left|f(x)-\widetilde H(x)\right|dx\leq
$$
$$
\leq\sum_{j=1}^{\nu_0}\int_{E_q^{(j)}}\left|\gamma_j\chi_{\Delta_j}(x)-\widetilde H_q^{(j)}(x)\right|dx+
$$
$$
+\sum_{j=1}^{\nu_0}\left(\sum_{l\in[1,\nu_0]\setminus\{j\}}\int_{[0,1]\setminus\Delta_l}\left|\widetilde H_q^{(l)}(x)\right|dx\right)<\frac{\varepsilon}{2}.
$$

Further, let $M$ be a natural number from $[2^{n_0},2^{n_{\nu_0}})$ . Then $M\in [2^{n_{m-1}},2^{n_m})$ for some $m\in[1,\nu_0]$. Taking (28)--(34) into account for any measurable set $e\subseteq E_\varepsilon$ we get
$$
\int_{e}\left| \sum_{k=2^{n_0}}^M \delta_k\tilde a_k W_k(x)\right|dx\leq
$$
$$
\leq\int_{e}\left| \sum_{j=1}^{m-1} \widetilde H_q^{(j)}(x)\right|dx+
\int_{e}\left| \sum_{k=2^{n_{m-1}}}^M \delta_k^{(m)}a_k^{(m)} W_k(x)\right|dx\leq
$$
$$
\leq \sum_{j=1}^{m-1}\int_{e}\bigl| \gamma_j\chi_{\Delta_j}(x)-\widetilde H_q^{(j)}(x)\bigr|dx+\int_{e}\left| \sum_{j=1}^{m-1}\gamma_j\chi_{\Delta_j}(x)\right|dx+
$$ 
$$
+\int_0^1\left| \sum_{k=2^{n_{m-1}}}^M \delta_k^{(m)}a_k^{(m)} W_k(x)\right|dx<
$$
$$
<\sum_{j=1}^{\nu_0}\sum_{l=1}^{\nu_0}\int_{E_q^{(l)}}\left|\gamma_j\chi_{\Delta_j}(x)-\widetilde H_q^{(j)}(x)\right|dx+\int_e\sum_{j=1}^{\nu_0}\bigl| \gamma_j\chi_{\Delta_j}(x)\big|dx+\frac{\varepsilon}{4}\leq
$$
$$
\leq\sum_{j=1}^{\nu_0}\int_{E_q^{(j)}}\left|\gamma_j\chi_{\Delta_j}(x)-\widetilde H_q^{(j)}(x)\right|dx
+\sum_{j=1}^{\nu_0}\left(\sum_{l\in[1,\nu_0]\setminus\{j\}}\int_{[0,1]\setminus\Delta_j}\left|\widetilde H_q^{(j)}(x)\right|dx\right)+
$$
$$
+\int_e|f(x)|dx+\frac{\varepsilon}{4}<\int_e|f(x)|dx+\frac{\varepsilon}{2}
$$
and
$$
\int_0^1\left| \sum_{k=2^{n_0}}^M \tilde a_k W_k(x)\right|dx\leq
$$
$$
\leq\sum_{j=1}^{\nu_0}\max_{2^{n_{j-1}}\leq N < 2^{n_j}}\int_0^1\left| \sum_{k=2^{n_{j-1}}}^N\tilde a_k^{(j)} W_k(x)\right|dx<\frac{\varepsilon}{2}.
$$
Hence, polynomials $\widetilde P(x)$ and $\widetilde H(x)$ satisfy all lemma 3 statements except for 1). To have strict inequalities between coefficients we choose a natural number $N_0>\log_2\frac{2}{\varepsilon}$ and set
$$
a_k = \tilde a_k + 2^{-(N_0+k)}.
$$
Now it is not hard to verify that new polynomials
$$
P(x)=\sum_{k=2^{n_0}}^{2^{n_{\nu_0}}-1}a_kW_k(x)\quad\hbox{and}\quad H(x)=\sum_{k=2^{n_0}}^{2^{n_{\nu_0}}-1}\delta_ka_kW_k(x),\quad \delta_k=\pm 1
$$
satisfy all lemma 3 statements. 

Lemma 3 is proved.

\paragraph{\bf Lemma 4} 
{ \it For any $\delta\in(0,1)$ there exists a weight function $0<\mu(x)\leq1$ with $|\{x\in[0,1];\ \mu(x)=1\}|>1-\delta$, so that for any number $n_0\in\mathbb{N}$, $\varepsilon\in(0,1)$ and step function $f(x)=\sum_{j=1}^{\nu_0} \gamma_j\chi_{\Delta_j}(x)$ with $\gamma_j\neq 0$ and $\sum_{j=1}^{\nu_0}|\Delta_j|=1$, where $\{\Delta_j\}_{j=1}^{\nu_0}$ are disjoint dyadic intervals, one can find polynomials 
$$
P(x)=\sum_{k=2^{n_{0}}}^{2^{n}-1}a_{k}W_{k}(x)\quad\hbox{and}\quad
H(x)=\sum_{k=2^{n_{0}}}^{2^{n}-1}\delta_{k}a_{k}W_{k}(x),\quad\delta_{k}=\pm1,
$$
in the Walsh system, which satisfy the following conditions:}
$$
0<a_{k+1}<a_{k}<\varepsilon ,\quad k\in\mathit{[\mathit{2}^{n_{0}},2^{n}-1)},\leqno1)
$$
$$
\int_0^1|f(x)-H(x)|\mu(x)dx<\varepsilon, \leqno2)
$$
$$
\max_{2^{n_{0}}\leq M<2^{n}}\int_0^1\left|\sum_{k=2^{n_{0}}}^{M}\delta_{k}
a_{k}W_{k}(x)\right|\mu(x)dx<\int_0^1|f(x)|\mu(x)dx+\varepsilon,\leqno3)
$$
$$
\max_{2^{n_{0}}\leq M<2^{n}}\int_0^1\left|\sum_{k=2^{n_{0}}}^{M}a_{k}
W_{k}(x)\right|dx<\varepsilon.\leqno4)
$$

\paragraph{\bf Proof of Lemma 4}
Let $\delta\in(0,1)$, $N_0 = 1$ and $\{f_m(x)\}_{m=1}^{\infty}$,
$$
f_m(x)=\sum_{j=1}^{\nu_m}\gamma_j^{(m)}\chi_{\Delta_j^{(m)}}(x),\leqno(35)
$$
be a sequence of all step functions with rational $\gamma_j^{(m)}\neq 0$ and $\sum_{j=1}^{\nu_m}|\Delta_j^{(m)}|=1$, where $\{\Delta_j^{(m)}\}_{j=1}^{\nu_m}$ are disjoint dyadic intervals.
By successively applying lemma 3 one can find sets $E_m\subset [0,1]$ and polynomials 
$$
P_{m}(x)=\sum_{k={2^{N_{m-1}}}}^{2^{N_m}-1}a_{k}^{(m)}W_{k}
(x),\leqno(36)
$$
$$
H_{m}(x)=\sum_{k={2^{N_{m-1}}}}^{2^{N_m}-1}\delta_{k}^{(m)}a_{k}
^{(m)}W_{k}(x),\quad \delta_{k}^{(m)}=\pm 1,\leqno(37)
$$
in the Walsh system which satisfy the following conditions for any natural number $m$:
$$
|E_m|>1-\frac{1}{2^{m+1}},\leqno(38)
$$
$$
0<a_{k+1}^{(m)}< a_k^{(m)}<\frac{1}{4^{N_{m-1}}}, \quad k\in[2^{N_{m-1}}, 2^{N_m}-1),\leqno(39)
$$ 
$$
\int_{E_m}\left|f_{m}(x)-H_{m}(x)\right|dx<\frac{1}{2^{m+1}},\leqno(40)
$$
$$
\max_{2^{N_{m-1}}\leq M<{2^{N_m}}}\int_e\left|\sum_{k={2^{N_{m-1}}}}^{M}
\delta_{k}^{(m)}a_{k}^{(m)}W_{k}(x)\right|dx<\int_e|f_m(x)|dx
+\frac{1}{2^{m+1}}\leqno(41)
$$
for any measurable set $e\subseteq E_m$ and
$$
\max_{2^{N_{m-1}}\leq M<2^{N_m}}\int_0^1\left|\sum_{k={2^{N_{m-1}}}}^{M}a_{k}^{(m)}W_{k}(x)\right|dx<\frac{1}{2^{m+1}}.\leqno(42)
$$

We set
$$
\begin{cases}
\Omega_{n}={\displaystyle\bigcap_{m=n}^{+\infty}}E_{m}
,\quad n\in\mathbb{N},
\\
E=\Omega_{\widetilde n}={\displaystyle\bigcap_{m=\widetilde n
}^{+\infty}}E_m,\quad \widetilde n=[\log_{1/2}\delta]+1,
\\
B=\Omega_{\widetilde n
}\bigcup\left(\displaystyle{\bigcup_{n=\widetilde n+1}^{+\infty}}\Omega_{n}
\setminus\Omega_{n-1}\right).
\end{cases}
\leqno(43)
$$

It is clear (see (38) and (43)) that
$$
|B|=1\quad\hbox{and}\quad |E|>1-\delta.
$$

We define a function $\mu(x)$ in the following way:
$$
\mu(x)=\begin{cases}
1,\quad x\in E\cup([0,1]\setminus B),
\\
\mu_{n},\quad \ x\in\Omega_{n}
\setminus\Omega_{n-1},\ \ n\geq \widetilde n+1,
\end{cases}
\leqno(44)
$$
where
$$
\mu_{n}=\frac{1}{2^{n+1}}\cdot\left[\prod_{m=1}^{n}
h_{m}\right]  ^{-1},\leqno(45)
$$

$$
h_m=1+\int_{0}^{1}
|f_m(x)|dx+\max_{2^{N_{m-1}}\leq M<2^{N_m}}\int_{0}^{1}\left|
\sum_{k=2^{N_{m-1}}}^{M}\delta_k^{(m)}a_{k}^{(m)}W_k(x)\right|dx.
$$

It follows  from (43)--(45) that for all $m\geq \widetilde n$
$$
\int_{[0,1]\setminus\Omega_{m}}\left|H_m(x)\right|
\mu(x)dx=\sum_{n=m+1}^{+\infty}\left(  \int_{\Omega_{n}\setminus\Omega_{n-1}
}\left|H_m(x)\right|\mu_{n}dx\right)<\leqno(46)
$$
$$
<\sum_{n=m+1}^{\infty}\frac{1}{2^{n+1}h_m}\left(  \int_{0}^{1}\left|H_{m}(x)\right|dx\right)<\frac{1}{2^{m+1}}.
$$

In a similar way for all $m\geq \widetilde n$ and $M\in[2^{N_{m-1}},2^{N_m})$ we have
$$
\int_{[0,1]\setminus\Omega_{m}}\left| f_{m}(x)\right|
\mu(x)dx<\frac{1}{2^{m+1}}\leqno(47)
$$
and
$$\int_{[0,1]\setminus\Omega_{m}}\left|\sum_{k=2^{N_{m-1}}}^{M}\delta_k^{(m)}a_k^{(m)}W_k(x)\right| \mu(x)dx< \frac{1}{2^{m+1}}.\leqno(48)$$

Since $\Omega_m\subset E_m$, then considering (40), (43)--(47), for all $m\geq \widetilde n$ we get
$$
\int_{0}^{1}|f_{m}(x)-{H}_{m}(x)|\mu(x)dx=\int_{\Omega_{m}}|f_{m}(x)-{H}_{m}(x)|\mu(x)dx+\leqno(49)
$$
$$
+\int_{[0,1]\setminus\Omega_{m}}|f_{m}(x)-{H}_{m}(x)|\mu
(x)dx<\frac{1}{2^{m+1}}+2\cdot\frac{1}{2^{m+1}}<
\frac{1}{2^{m-1}}.
$$

Further, taking relations (41), (43)--(45), (48) into account for all  $M\in[2^{N_{m-1}},2^{N_m})$ and $m\geq \widetilde n+1$ we have
$$
\int_{0}^{1}\left|\sum_{k=2^{N_{m-1}}}^{M}\delta_k^{(m)}a_k^{(m)}W_k(x)\right| \mu(x)dx=\leqno(50)
$$
$$
=\int_{\Omega_{m}}\left|\sum_{k=2^{N_{m-1}}}^{M}\delta_k^{(m)}a_k^{(m)}W_k(x)\right|\mu(x)dx+
$$
$$
+\int_{[0,1]\setminus\Omega_{m}}\left|\sum_{k=2^{N_{m-1}}}^{M}\delta_k^{(m)}a_k^{(m)}W_k(x)\right| \mu(x)dx<
$$
$$
<\int_{\Omega_{\widetilde n}}\left|\sum_{k=2^{N_{m-1}}}^{M}\delta_k^{(m)}a_k^{(m)}W_k(x)\right|\mu(x)dx+
$$
$$
+\sum_{n=\widetilde n+1}^{m}\mu_{n}\cdot\int_{\Omega_{n}\setminus\Omega_{n-1}}\left|
\sum_{k=2^{N_{m-1}}}^{M}\delta_k^{(m)}a_k^{(m)}W_k(x)\right|dx+\frac{1}{2^{m+1}}<
$$
$$
<\int_{\Omega_{\widetilde n}}|f_m(x)|dx
+\frac{1}{2^{m+1}}+\sum_{n=\widetilde n+1}^{m}\mu_{n}\left(\int_{\Omega_{n}\setminus\Omega_{n-1}}|f_m(x)|dx+\frac{1}{2^{m+1}}\right)+\frac{1}{2^{m+1}}=
$$
$$
=\int_{\Omega_{\widetilde n}}|f_{m}(x)|dx+\sum_{n=\widetilde n+1}^{m}\int_{\Omega_{n} \setminus\Omega_{n-1}}|f_{m}(x)|\cdot\mu_{n}dx+
$$
$$
+\frac{1}{2^{m+1}}\left(2+\sum_{n=\widetilde n+1}^{m}\mu_n\right)<\int_0^1|f_m(x)|\mu(x)dx+\frac{1}{2^{m-1}}.
$$

Let $n_0\in\mathbb{N}$ and $\varepsilon\in(0,1)$ be arbitrarily given. From the sequence (35) we choose such a function $f_{m_0}(x)$ that 
$$
m_0>\max\left\{\tilde n,\ \log_2\frac{8}{\varepsilon}\right\}\quad\hbox{and}\quad 2^{N_{m_0-1}}>2^{n_0}, \leqno(51)
$$
$$
\int_0^1|f(x)-f_{m_{0}}(x)|dx<\frac{\varepsilon}{4}\leqno(52)
$$
and for $k\in\left[2^{n_0},2^{N_{m_0}}\right)$ set
$$
a_k=
\begin{cases}
a_{2^{N_{m_0-1}}}^{(m_0)}+\frac{1}{2^{k+m_0}},& \hbox{when} \quad k\in\left[2^{n_0},2^{N_{m_0-1}}\right),
\\
\\
a_k^{(m_0)},& \hbox{when} \quad k\in\left[2^{N_{m_0-1}},2^{N_{m_0}}\right),
\end{cases}\leqno(53)
$$
$$
\delta_k=\begin{cases}
1,& \hbox{when} \quad k\in\left[2^{n_0},2^{N_{m_0-1}}\right),
\\
\\
\delta_k^{(m_0)}=\pm 1,& \hbox{when} \quad k\in\left[2^{N_{m_0-1}},2^{N_{m_0}}\right),
\end{cases}\leqno(54)
$$
and
$$
P(x)=\sum_{k=2^{n_0}}^{2^{N_{m_0}}-1}a_kW_k(x)=\sum_{k=2^{n_0}}^{2^{N_{m_0-1}}-1}a_kW_k(x)+P_{m_0}(x),
$$
$$
H(x)=\sum_{k=2^{n_0}}^{2^{N_{m_0}}-1}\delta_ka_kW_k(x)=\sum_{k=2^{n_0}}^{2^{N_{m_0-1}}-1}a_kW_k(x)+H_{m_0}(x).
$$

Now let us verify that the function $\mu(x)$ and polynomials $P(x)$ and $H(x)$ satisfy all requirements of lemma 4.
The statement 1) immediately follows from (39), (51) and (53). By referring to (1.a), (39), (49), (51)--(54) we get the statement 2):
$$
\int_{0}^{1}|f(x)-{H}(x)|\mu(x)dx\leq\int_{0}^{1}|f(x)-f_{m_0}(x)|dx+
$$
$$
+\int_{0}^{1}|f_{m_0}(x)-H_{m_0}(x)|\mu(x)dx+a_{2^{N_{m_0-1}}}^{(m_0)}\cdot\int_0^1\left|\sum_{k=2^{n_0}}^{2^{N_{m_0-1}}-1}W_k(x)\right|dx+
$$
$$
+\sum_{k=2^{n_0}}^{2^{N_{m_0-1}}-1}\frac{1}{2^{k+m_0}}<\varepsilon.
$$

Further, by using  (50)--(54) we obtain
$$
\max_{2^{n_0}\leq M<{2^{N_{m_0}}}}\int_0^1\left|\sum_{k=2^{n_0}}^M\delta_ka_kW_k(x)\right|\mu(x)dx\leq
$$
$$
\leq\max_{2^{n_0}\leq M_1<{2^{N_{m_0-1}}}}\int_0^1\left|\sum_{k=2^{n_0}}^{M_1}a_kW_k(x)\right|dx+
$$ 
$$
+\max_{{2^{N_{m_0-1}}}\leq M_2<{2^{N_{m_0}}}}\int_0^1\left|\sum_{k={2^{N_{m_0-1}}}}^{M_2}\delta_k^{(m_0)}a_k^{(m_0)}W_k(x)\right|\mu(x)dx<
$$
$$
<\max_{2^{n_0}\leq M_1<{2^{N_{m_0-1}}}}\int_0^1\left|\sum_{k=2^{n_0}}^{M_1}a_kW_k(x)\right|dx+\int_0^1|f(x)|\mu(x)dx+\frac{\varepsilon}{2}.
$$
Let $M_1$ be an arbitrary natural number from $\left[2^{n_0},2^{N_{m_0-1}}\right)$. Then $M_1\in [2^{n_1},2^{n_1+1})$ for some $n_1\in[n_0,N_{m_0-1})$ and, considering (1.a), we have 
$$
\int_0^1\left|\sum_{k=2^{n_0}}^{M_1}a_kW_k(x)\right|dx< a_{2^{N_{m_0-1}}}^{(m_0)}\cdot\int_0^1\left|\sum_{k=2^{n_0}}^{2^{n_1}-1}W_k(x)\right|dx+a_{2^{N_{m_0-1}}}^{(m_0)}\cdot2^{n_1}+
$$
$$
+\sum_{k=2^{n_0}}^{M_1}\frac{1}{2^{k+m_0}}<\frac{\varepsilon}{2},
$$
which proves the statement 3). The statment 4) can be proved in a similar way by using relations (42), (51) and (53). 

Lemma 4 is proved.

\section{Proof of  theorem}
\label{Thm}

Let $\delta\in(0,1)$ and $\{f_m(x)\}_{m=1}^{\infty}$,
$$
f_m(x)=\sum_{j=1}^{\nu_m}\gamma_j^{(m)}\chi_{\Delta_j^{(m)}}(x),\leqno(55)
$$
be a sequence of all step functions with rational $\gamma_j^{(m)}\neq 0$ and $\sum_{j=1}^{\nu_m}|\Delta_j^{(m)}|=1$, where $\{\Delta_j^{(m)}\}_{j=1}^{\nu_m}$ are disjoint dyadic intervals. 

By applying lemma 4 we obtain a weight function $0<\mu(x)\leq1$ with $|\{x\in[0,1],\ \mu(x)=1\}|>1-\delta$ and polynomials 
$$
P_{m}(x)=\sum_{k={N_{m-1}}}^{{N_{m}^{{}}}-1}a_{k}^{(m)}W_{k}
(x),\leqno(56)
$$
$$
H_{m}(x)=\sum_{k={N_{m-1}}}^{{N_{m}^{{}}}-1}\delta_{k}^{(m)}a_{k}%
^{(m)}W_{k}(x),\quad \delta_{k}^{(m)}=\pm 1,\leqno(57)
$$
in the Walsh system, which satisfy the following conditions for each natural number $m$:
$$
\begin{cases}
0<a_{k+1}^{(1)}<a_{k}^{(1)},
\\
0<a_{k+1}^{(m)}<a_{k}^{(m)}<\min\bigl\{2^{-m}, a_{{N_{m - 1}} - 1}^{(m-1)}\bigr\}, & \hbox{for}\quad m>1,
\end{cases}\leqno(58)
$$
when $ k\in[N_{m-1},N_m-1)$,
$$
\int_0^1\left|f_{m}(x)-H_{m}(x)\right|\mu(x)dx<2^{-m-2},\leqno(59)
$$
and for each natural number $M\in[N_{m-1},{N_m})$
$$
\int_0^1\left|\sum_{k={N_{m-1}}}^{M}
\delta_{k}^{(m)}a_{k}^{(m)}W_{k}(x)\right|\mu(x)dx<\int_0^1|f_m(x)|\mu(x)dx
+2^{-m},\leqno(60)
$$
and
$$
\int_0^1\left|\sum_{k=N_{m-1}}^{M}a_{k}^{(m)}W_{k}(x)\right|dx<2^{-m-2}.\leqno(61)
$$

From (56) and (61) it immediately follows that
$$
\int_{0}^{1}\left|\sum_{m=1}^{\infty}P_{m}(x)\right|dx\leq\sum_{m=1}^{\infty}\int_{0}^{1}|P_{m}(x)|dx<+\infty.\leqno(62)
$$

We denote $P_0(x)=\sum_{k=0}^{N_0-1}a_{k}^{(0)}W_{k}(x)$,
where coefficients $a_k^{(0)}$ are arbitrary monotonically decreasing positive numbers with $a_{N_0-1}^{(0)}>a_{N_0}^{(1)}$ and define a function $g(x)$ and a series $\sum_{k=0}^{\infty}a_kW_k(x)$ in the following way:
$$
g(x)=\sum_{m=0}^{\infty}P_{m}(x),\leqno(63)
$$
$$
a_{k}=a_{k}^{(m)},\ \mbox{when} \ k\in[N_{m-1}, N_m),\ m\geq 0,\ N_{-1}=0.\leqno(64)
$$

Considering (58), (61)--(64) we conclude that the series $\sum_{k=0}^{\infty}a_kW_k(x)$ converges to $g\in L^1[0,1]$ in $L^1[0,1]$ metric, and $a_{k}=\int_{0}^{1}g(x)W_{k}(x)dx \searrow0$.

Further, let $f(x)$ be an arbitrary function from $L^1_\mu[0,1]$. Then one can choose such a polynomial $f_{\nu_{1}}(x)$ from the sequence (55), that that
$$
\int_0^1\left|f(x)-P_0(x)-f_{\nu_{1}}(x)\right|\mu(x)dx<2^{-3}.\leqno(65)
$$

By denoting  
$$\delta_{k}= 
\begin{cases}
\delta_{k}^{({\nu_{1}})}=\pm 1,&  \hbox{when}\quad k\in[N_{{\nu}_1-1},N_{\nu_1})
\\
1, &  \hbox{when}\quad k\in[0,N_{{\nu}_1-1})
\end{cases}$$
and taking (57), (59), (60) and (65) into account we have
$$
\int_0^1\left|f(x)-\sum_{k=0}^{N_{\nu_1}-1}\delta_{k}a_{k}
W_{k}(x)\right|\mu(x)dx=
$$
$$
=\int_0^1\left|f(x)-\sum_{m=0}^{\nu_1-1}P_m(x) - H_{\nu_1}(x)\right|\mu(x)dx<
$$
$$
<\int_0^1\left|f(x)-P_0(x)-f_{\nu_{1}}(x)\right|\mu(x)dx+\int_0^1\left|f_{\nu_1}(x)-H_{\nu_1}(x)\right|\mu(x)dx
+
$$
$$
+\sum_{m=1}^{\nu_1-1}\int_{0}^{1}|P_m(x)|dx<2^{-3}+2^{-\nu_1-2}+\sum_{m=1}^{\nu_1-1}2^{-m-2}<2^{-1},
$$
and for any natural number $M\in[N_{{\nu}_1-1},N_{{\nu}_1})$ we get
$$
\int_0^1\left|\sum_{k=N_{\nu_1-1}}^{M}\delta_{k}a_{k}
W_{k}(x)\right|\mu(x)dx<\int_0^1|f_{\nu_{1}}(x)|\mu(x)dx+2^{-\nu_1-2}.
$$

Assume that for $q>1$ numbers $\nu_{1}<\nu_{2}<\dots<\nu_{q-1}$ and $\{\delta_k=\pm 1\}_{k=0}^{N_{\nu_{q-1}}-1}$ are already chosen, so that for each natural number $j\in[1,q-1]$ the following conditions hold:
$$\delta_{k}= 
\begin{cases}
\delta_{k}^{({\nu_{j}})}=\pm 1,&  \hbox{when}\quad k\in[N_{{\nu}_j-1},N_{\nu_j}),
\\
1, &  \hbox{when}\quad  k\notin\bigcup_{j=1}^{q-1}[N_{{\nu}_j-1},N_{\nu_j}),
\end{cases}
$$
$$
\int_0^1\left|f(x)-\sum_{k=0}^{N_{\nu_j}-1}\delta_{k}a_{k}
W_{k}(x)\right|\mu(x)dx<2^{-j},\leqno(66)
$$
and for any natural number $M\in[N_{{\nu}_j-1},N_{{\nu}_j})$
$$
\int_0^1\left|\sum_{k=N_{\nu_j-1}}^{M}\delta_ka_k
W_{k}(x)\right|\mu(x)dx<\int_0^1|f_{\nu_{j}}(x)|\mu(x)dx+2^{-\nu_j-2}.
$$

We choose a function $f_{\nu_{q}}(x)$ from the sequence (55) with $\nu_{q}>\nu_{q-1}$ so that
$$
\int_0^1\left|\left[f(x)-\sum_{k=0}^{N_{\nu_{q-1}}-1}\delta_{k}a_{k}
W_{k}(x)\right]-f_{\nu_q}(x)\right|\mu(x)dx<2^{-q-2},\leqno(67)
$$
and define
$$\delta_{k}= 
\begin{cases}
\delta_{k}^{({\nu_{q}})}=\pm 1,&  \hbox{when}\quad k\in\bigr[N_{{\nu}_q-1},N_{\nu_q}\bigr),
\\
1, &  \hbox{when}\quad  k\notin\bigcup_{j=1}^{q}[N_{{\nu}_j-1},N_{\nu_j}).
\end{cases}\leqno(68)$$ 

Taking (59), (61), (64), (67) and (68) into account we obtain
$$
\int_0^1\left|f(x)-\sum_{k=0}^{N_{\nu_q}-1}\delta_{k}a_{k}
W_{k}(x)\right|\mu(x)dx=\leqno(69)
$$
$$
=\int_0^1\left|f(x)-\sum_{k=0}^{N_{\nu_{q-1}}-1}\delta_{k}a_{k}
W_{k}(x)-\sum_{m=\nu_{q-1}+1}^{\nu_q-1}P_m(x) - H_{\nu_q}(x)\right|\mu(x)dx<
$$
$$
<\int_0^1\left|\left[f(x)-\sum_{k=0}^{N_{\nu_{q-1}}-1}\delta_{k}a_{k}
W_{k}(x)\right]-f_{\nu_q}(x)\right|\mu(x)dx+
$$
$$
+\int_0^1\left|f_{\nu_q}(x)-H_{\nu_q}(x)\right|\mu(x)dx
+\sum_{m=\nu_{q-1}+1}^{\nu_q-1}\int_{0}^{1}|P_m(x)|dx\leq
$$
$$
\leq2^{-q-2}+2^{-\nu_q-2}+\sum_{m=\nu_{q-1}+1}^{\nu_q-1}2^{-m-2}< 2^{-q-2}+2^{-q-2}+2^{-q-1}=2^{-q}.
$$

Further, from (66) and (67) we have
$$
\int_0^1|f_{\nu_{q}}(x)|\mu(x)dx<\int_0^1\left|\left[f(x)-\sum_{k=0}^{N_{\nu_{q-1}}-1}\delta_{k}a_{k}
W_{k}(x)\right]-f_{\nu_q}(x)\right|\mu(x)dx+
$$
$$
+\int_0^1\left|f(x)-\sum_{k=0}^{N_{\nu_{q-1}}-1}\delta_{k}a_{k}
W_{k}(x)\right|\mu(x)dx<2^{-q-2} + 2^{-q+1}<2^{-q+2}.
$$
Thus, from (60) and (68) it follows that for each natural number $M\in[N_{{\nu}_q-1},N_{{\nu}_q})$
$$
\int_0^1\left|\sum_{k=N_{{\nu}_q-1}}^{M}\delta_ka_k
W_{k}(x)\right|\mu(x)dx<\int_0^1|f_{\nu_{q}}(x)|\mu(x)dx+2^{-\nu_q-2}<2^{-q+3}\leqno(70)
$$
is true.

Apparently, by using induction one can determine a growing sequence of natural numbers $\{\nu_q\}_{q=1}^\infty$ and numbers $\delta_k=\pm 1$ so that the conditions (68)--(70) hold for any $q\in\mathbb{N}$. Hence, considering also (61), we obtain a series
$$
\sum_{k=0}^{\infty}\delta_{k}a_{k}W_{k}(x),\quad \delta_k=\pm 1,
$$
which converges to $f$ in $L^1_\mu[0,1]$ metric. 

The theorem is proved.

\end{document}